\documentclass{article}

\usepackage{arxiv}

\usepackage[utf8]{inputenc} 
\usepackage[T1]{fontenc}    
\usepackage{hyperref}       
\usepackage{url}            
\usepackage{booktabs}       
\usepackage{amsfonts}       
\usepackage{nicefrac}       
\usepackage{microtype}      
\usepackage{lipsum}		
\usepackage{graphicx}
\usepackage{natbib}
\usepackage{doi}
\usepackage{bm}
\usepackage{graphicx}%
\usepackage{multirow}%
\usepackage[title]{appendix}%
\usepackage{graphicx}%
\usepackage{multirow}%
\usepackage{amsmath,amssymb,amsfonts}%
\usepackage{amsthm}%
\usepackage{mathrsfs}%
\usepackage[title]{appendix}%
\usepackage{xcolor}%
\usepackage{textcomp}%
\usepackage{manyfoot}%
\usepackage{booktabs}%
\usepackage{algorithm}%
\usepackage{algorithmicx}%
\usepackage{algpseudocode}%
\usepackage{listings}%
\usepackage{bm}
\usepackage{etoolbox}
\usepackage{graphicx}
\usepackage{array}
\apptocmd{\thebibliography}{\csname phantomsection\endcsname}{}{}
\usepackage{graphicx}
\usepackage{subcaption}
\usepackage{float}
\usepackage{multirow}
\usepackage{setspace}
\usepackage{geometry}

\title{Discontinuous Behavior of Time-of-Flight Distribution for Bi-impulsive Earth-Moon Transfers in the Three-Body Model}

\author{
 Shuyue Fu \\
  School of Astronautics\\
  Beihang University\\
  Beijing, 100191 \\
  \texttt{fushuyue@buaa.edu.cn} \\
   \And
  Di Wu \\
  School of Astronautics\\
  Beihang University\\
  Beijing, 100191 \\
  \texttt{wudi2025@buaa.edu.cn} \\
  \And
  Shengping Gong \\
  School of Astronautics\\
  Beihang University\\
  Beijing, 100191 \\
  \texttt{gongsp@buaa.edu.cn} \\
}

\begin{document}
\maketitle
\begin{abstract}
As interest in the Earth-Moon transfers renewed around the world, understanding the solution space of transfer trajectories facilitates the construction of transfers. This paper is devoted to reporting a novel or less-reported phenomenon about the solution space of bi-impulsive Earth-Moon transfers in the Earth-Moon planar circular restricted three-body problem. Differing from the previous works focusing on the transfer characteristics of the solution space, we focus on the distribution of the construction parameters, i.e., departure phase angle at the Earth parking orbit, initial-to-circular velocity ratio, and time of flight. Firstly, the construction method of bi-impulsive transfers is described, and the solutions satisfying the given constraints are obtained from the grid search method and trajectory correction. Then, the distribution of the obtained solutions is analyzed, and an interesting phenomenon about the discontinuous behavior of the time-of-flight distribution for each departure phase angle is observed and briefly reported. This phenomenon can further provide useful insight into the construction of bi-impulsive transfers, deepening the understanding of the corresponding solution space. 
\end{abstract}


\section{Introduction}\label{sec1}
Recently, the interest in the Earth-Moon transfers has been renewed due to the proposal of several lunar ex-ploration missions \cite{zheng2023analysis}. In the preliminary design stage, the Earth-Moon transfer (with multiple impulses) can be simplified into a bi-impulsive scenario where the spacecraft performs an Earth injection impulse to de-part from the Earth parking orbit and inserts into the Moon target orbit by a Moon insertion impulse \cite{topputo2013optimal,grossi2024optimal}. The classical scenario is the bi-impulsive scenario from a circular Earth parking orbit to a circular Moon target orbit \cite{topputo2013optimal}. To construct transfers in this scenario, the first step is to select the dynamical models. Compared to direct transfers in the patched two-body problem \cite{battin1999introduction}, there is a potential to use the multi-body problems (e.g., the Earth-Moon planar three-body problem (PCR3BP) \cite{fu2025escape}, the Sun-Earth/Moon planar bicircular restricted four-body problem (PBCR4BP) \cite{topputo2013optimal,grossi2024optimal}, or other high-er-fidelity models \cite{wang2024low}) to obtain more solutions, including direct and low-energy transfers, to satisfy specific mission requirements. Therefore, we adopt one of the multi-body problems, i.e., the Earth-Moon PCR3BP, to construct the bi-impulsive Earth-Moon transfers and analyze the corresponding distribution.

Because there is no closed-form solution in the Earth-Moon PCR3BP, numerical methods are adopted to construct bi-impulsive transfers. Typically, the construction method can be divided into two steps: generation of initial guesses and trajectory correction. Firstly, the initial guesses of transfer trajectories are generated by several methods, such as grid search \cite{topputo2013optimal,oshima2019low} and methods based on the dynamical system theory \cite{scheuerle2025energy}. Then, initial guesses are corrected to satisfy the specific constraints, i.e., the orbital altitude of the Earth parking orbit/Moon target orbit and tangential departure/insertion \cite{topputo2013optimal}. Compared to the methods based on the dynamical system theory, the grid search method typically yields more solutions and provides a global map to analyze the transfer characteristics \cite{topputo2013optimal,oshima2019low}. Topputo \cite{topputo2013optimal} and Oshima et al. \cite{oshima2019low} presented a global map of impulse and time of flight (TOF) for the bi-impulsive transfers within 100 and 200 Days in the Sun-Earth/Moon PBCR4BP, revealing the rich structure of the solution space. Based on their work, Campana and Topputo \cite{campana2024clustering} used the data mining method to extract the transfer families in the Sun-Earth/Moon PBCR4BP. The aforementioned works \cite{topputo2013optimal,oshima2019low,campana2024clustering} performed a valuable exploration on the solution space and trans-fer characteristics of the bi-impulsive Earth-Moon transfers in the multi-body problems. However, they did not focus on the distribution of the construction parameters (i.e., parameters determining the transfer trajectories, such as the departure phase angle at the Earth parking orbit, initial-to-circular velocity ratio, and TOF \cite{topputo2013optimal}). The analysis of these distributions can further provide useful insight into selecting suitable initial guesses (generated by construction parameters) to aid in the construction of bi-impulsive transfers. Therefore, in this paper, we adopt the grid search method and trajectory correction to obtain transfers in the Earth-Moon PCR3BP, and further explore the cor-responding distribution of the construction parameters.
Firstly, the construction method of bi-impulsive transfers in the Earth-Moon PCR3BP is described, and the obtained solutions are presented. Then, we focus on the distribution of the construction parameters, and report an interesting, novel, or less-reported phenome-non about the discontinuous behavior of TOF distribution for each departure phase angle. Although the map of TOF and departure phase angle in the Sun-Earth/Moon PBCR4BP was presented by Pinelli \cite{pinelli2023neural}, the global solution space was not investigated. Therefore, the discontinuous behavior of the TOF distribution was not clearly revealed. This paper indicates that such behavior can be clearly observed in the Earth-Moon PCR3BP, while in the Sun-Earth/Moon PBCR4BP, this type of distribution is expected to be more scattered due to the additional construction pa-rameter, i.e., the initial solar phase angle \cite{pinelli2023neural}. Based on the aforementioned discussion, the main contribution of this paper is to report this discontinuous behavior of TOF distribution for bi-impulsive transfers in the Earth-Moon PCR3BP, which can be considered as a novel or less-reported phenomenon. This phenomenon can further provide useful insight into the construction of bi-impulsive Earth-Moon transfers in the Earth-Moon PCR3BP.

The rest of this paper is organized as follows. Section \ref{sec2} introduces the background of this paper. The construction method of bi-impulsive transfers is described in Section \ref{sec3}. The obtained solutions and the novel (or less-reported) phenomenon about discontinuous behavior of TOF distribution for each departure phase angle are presented in Section \ref{sec4}. Finally, Section \ref{sec5} presents the conclusions of this paper.

\section{Mathematical Background}\label{sec2}
\subsection{Earth-Moon PCR3BP}\label{subsec2.1}
In this paper, the Earth-Moon PCR3BP is adopted to investigate the bi-impulsive Earth-Moon transfers. In this model, the motion of the spacecraft is dominated by the gravity forces from both the Earth and the Moon. The Earth-Moon rotating frame and dimensionless units are adopted to describe the equation of motion. The mass unit (MU) is set as the combined mass of the Earth and the Moon, the length unit (LU) is set as the distance between the Earth and the Moon, and the time unit (TU) is set as
$T_{\text{EM}}/2\pi$, where $T_{\text{EM}}$ denotes the orbital period of the Earth and the Moon around their barycenter. The parameter setting of this model can be found in Ref. \cite{fu2025escape}. Then, the equation of motion can be expressed as:
\begin{equation}
\left[ {\begin{array}{*{20}{c}}
  {\begin{array}{*{20}{c}}
  {\dot x} \\ 
  {\dot y} 
\end{array}} \\ 
  {\begin{array}{*{20}{c}}
  {\dot u} \\ 
  {\dot v} 
\end{array}} 
\end{array}} \right] = \left[ {\begin{array}{*{20}{c}}
  {\begin{array}{*{20}{c}}
  u \\ 
  v 
\end{array}} \\ 
  {\begin{array}{*{20}{c}}
  {2v + \frac{{\partial {\Omega _3}}}{{\partial x}}} \\ 
  { - 2u + \frac{{\partial {\Omega _3}}}{{\partial y}}} 
\end{array}} 
\end{array}} \right]\label{eq1}
\end{equation}
\begin{equation}
{\Omega _3} = \frac{1}{2}\left( {{x^2} + {y^2}} \right) + \frac{{1 - \mu }}{{{r_1}}} + \frac{\mu }{{{r_2}}}\label{eq2}
\end{equation}
where $bm{X} = {\left[ {x,{\text{ }}y,{\text{ }}u,{\text{ }}v} \right]^{\text{T}}}$ is the orbital state in the Earth-Moon rotating frame, $\mu$ is the mass parameter calculated by $\mu  = {m_2}/\left( {{m_1} + {m_2}} \right)$ (where $m_1$ and $m_2$ denote the masses of the Earth and the Moon), and $\Omega_3$ is the effective potential of the Earth-Moon PCR3BP. Furthermore, the distances between the spacecraft and the Earth ($r_1$) and the Moon ($r_2$) are expressed as:
\begin{equation}
{r_1} = \sqrt {{{\left( {x + \mu } \right)}^2} + {y^2}}\label{eq3}
\end{equation}
\begin{equation}
{r_2} = \sqrt {{{\left( {x + \mu -1} \right)}^2} + {y^2}}\label{eq4}
\end{equation}

Since Eq. \eqref{eq1} has no closed-form solution, we adopt MATLAB®'s ode113 command with absolute and relative tolerances set to $1\times10^{-13}$ to integrate trajectories in the Earth-Moon PCR3BP. Subsequently, the concept of bi-impulsive Earth-Moon transfers is introduced.

\subsection{ Bi-impulsive Earth-Moon Transfer}\label{subsec2.2}
Bi-impulsive Earth-Moon transfer describes a scenario where the spacecraft departs from a circular Earth parking orbit (whose orbital altitude is denoted as $h_i$) after an Earth injection impulse ($\Delta v_i$), coasts in the transfer trajectory, and finally performs a Moon insertion impulse ($\Delta v_f$) to insert into a circular Moon target orbit (whose orbital altitude is denoted as $h_f$) \cite{topputo2013optimal,oshima2019low}. Similar to Ref. \cite{topputo2013optimal,oshima2019low}, we also adopt the tangential impulses to construct the transfer trajectories. Therefore, the constraints of the trajectories should satisfy can be expressed as:
\begin{equation}
{\bm{\psi }_i} = \left[ {\begin{array}{*{20}{c}}
  {{{\left( {{x_i} + \mu } \right)}^2} + {y_i}^2 - {{\left( {{R_{\text{E}}} + {h_i}} \right)}^2}} \\ 
  {\left( {{x_i} + \mu } \right)\left( {{u_i} - {y_i}} \right) + {y_i}\left( {{v_i} + {x_i} + \mu } \right)} 
\end{array}} \right] = \bm{0}
\label{eq5}
\end{equation}
\begin{equation}
{\bm{\psi }_f} = \left[ {\begin{array}{*{20}{c}}
  {{{\left( {{x_f} + \mu  - 1} \right)}^2} + {y_f}^2 - {{\left( {{R_{\text{M}}} + {h_f}} \right)}^2}} \\ 
  {\left( {{x_f} + \mu  - 1} \right)\left( {{u_f} - {y_f}} \right) + {y_f}\left( {{v_f} + {x_f} + \mu  - 1} \right)} 
\end{array}} \right] = \bm{0}
\label{eq6}
\end{equation}
where the subscripts “\textit{i}” and “\textit{f}” denote the quantities corresponding to the departure and insertion points at the Earth parking orbit and Moon target orbit, respectively. In this paper, following Refs. \cite{topputo2013optimal,oshima2019low}, we also investigate bi-impulsive transfers with $h_i=167\text{ km}$ and $h_f=100\text{ km}$. The radius of the Earth and the Moon are set as $R_\text{E}=6378.145 \text{ km}$ and $R_\text{M}=1737.1 \text{ km}$. When the solution satisfying the constraints Eqs. \eqref{eq5}-\eqref{eq6} is obtained, the impulses of transfer can be calculated by:
\begin{equation}
\Delta {v_i} = \sqrt {{{\left( {{u_i} - {y_i}} \right)}^2} + {{\left( {{v_i} + {x_i} + \mu } \right)}^2}}  - \sqrt {\frac{{1 - \mu }}{{{R_{\text{E}}} + {h_i}}}} 
\label{eq7}
\end{equation}
\begin{equation}
\Delta {v_f} = \sqrt {{{\left( {{u_f} - {y_f}} \right)}^2} + {{\left( {{v_f} + {x_f} + \mu  - 1} \right)}^2}}  - \sqrt {\frac{\mu }{{{R_{\text{M}}} + {h_f}}}}
\label{eq8}
\end{equation}
\begin{equation}
\Delta v = \Delta {v_i} + \Delta {v_f}
\label{eq9}
\end{equation}
where $\Delta v$ denotes the total impulse of the bi-impulsive transfer. After introducing the bi-impulsive transfer, we present the construction method of this type of transfer.

\section{Construction of Bi-impulsive Transfers}\label{sec3}
\subsection{Construction Parameter}\label{subsec3.1}
Similar to Ref. \cite{topputo2013optimal}, we set the construction parameters of bi-impulsive transfers as:
\begin{equation}
\bm{y} = {\left[ {{\alpha _i},{\text{ }}{\beta _i},{\text{ TOF}}} \right]^{\text{T}}}
\label{eq10}
\end{equation}
where $\alpha _i$ denotes the departure phase angle at the Earth parking orbit, $\beta _i$ denotes the initial-to-circular velocity ratio, and TOF denotes the time of flight (TOF). With this setting, the orbital state of the departure point can be expressed as (assuming that the spacecraft departs from the prograde Earth parking orbit):
\begin{equation}
\left\{ \begin{gathered}
  {x_i} = {r_i}\cos {\alpha _i} - \mu  \hfill \\
  {y_i} = {r_i}\sin {\alpha _i} \hfill \\
  {u_i} =  - \left( {{\beta _i}\sqrt {\frac{{1 - \mu }}{{{r_i}}}}  - {r_i}} \right)\sin {\alpha _i} \hfill \\
  {v_i} = \left( {{\beta _i}\sqrt {\frac{{1 - \mu }}{{{r_i}}}}  - {r_i}} \right)\cos {\alpha _i} \hfill \\ 
\end{gathered}  \right.
\label{eq11}
\end{equation}
where $r_i=R_\text{E}+h_i$. The orbital state calculated by Eq. \eqref{eq8} satisfies the constraint Eq. \eqref{eq5} rigorously. Then, we use this setting to generate the initial guesses of bi-impulsive transfers.
\subsection{Initial Guess}\label{subsec3.2}
In this paper, we adopt the grid search method to generate the initial guesses. We set $\bm{y}$ in Eq. \eqref{eq7} as ${\alpha _i} \in \left[ {0,{\text{ }}2\pi } \right){\text{ rad}}$ with a step-size of $\pi /36{\text{ rad}}$, ${\beta _i} \in \left[ {1.4,{\text{ }}1.414} \right]$ with a step-size of 0.0002, and ${\text{TOF}} \in \left[ {\pi /50,{\text{ 10}}\pi } \right]{\text{ TU}}$ with a step-size of $\pi /50{\text{ TU}}$. The ranges of $\alpha_i$ and $\beta_i$ are selected according to Ref. \cite{topputo2013optimal}, and the selection of TOF does not affect the qualitative findings reported in this paper. Therefore, the initial guess is generated from the initial state calculated by Eq. \eqref{eq8} and numerical integration within the time interval $\left[ {0,{\text{ TOF}}} \right]$.
\subsection{Trajectory Correction}\label{subsec3.3}
Once the initial guesses are generated, trajectory correction is performed to make the trajectories satisfy the constraints Eqs. \eqref{eq5}-\eqref{eq6}. Trajectory correction can be transformed into a nonlinear programming (NLP) problem. In this paper, we use MATLAB®'s fmincon com-mand to solve this NLP problem. To focus on the de-pendency of TOF on $\alpha_i$ for the obtained solutions, we fix the value of $\alpha_i$ during solving the NLP problem, i.e., $\alpha_i$ is not adjusted for each initial guess when performing the fmincon command. Therefore, the NLP variables are set as $\beta_i$ and TOF. Since the initial guesses generated by the method mentioned in Section \ref{subsec3.2} satisfy the constraint Eq. \eqref{eq5} rigorously, the objective function and constraint of the NLP problem are both set as Eq. \eqref{eq6}. When performing the fmincon command to solve this NLP problem, we use the sequential quadratic programming method, and adopt the following setting:
\begin{equation}
\left\{ \begin{gathered}
  {\beta _i}_{{\text{lower boundary}}} = 1.4 \hfill \\
  {\beta _i}_{{\text{upper boundary}}} = 1.414 \hfill \\
  {\text{TO}}{{\text{F}}_{{\text{lower boundary}}}} = \frac{\pi }{{50}} \hfill \\
  {\text{TO}}{{\text{F}}_{{\text{upper boundary}}}} = 10\pi  \hfill \\
  {\text{TolX}} = 1 \times {10^{ - 8}} \hfill \\
  {\text{TolFun}} = 1 \times {10^{ - 8}} \hfill \\
  {\text{TolCon}} = 1 \times {10^{ - 8}} \hfill \\
  {\text{MaxIter}} = 500{\text{    MaxFunEvals}} = 500 \hfill \\ 
\end{gathered}  \right.
\label{eq12}
\end{equation}

During solving the NLP problem, the Earth/Moon collision trajectories are excluded. When the correction is finished, the solutions satisfying $\left\| {{\bm{\psi }_f}} \right\| < 1 \times {10^{ - 8}}$ are selected and recorded. Subsequently, we present the obtained solutions and analyze the distribution of the corresponding construction parameters.

\section{Results and Discussion}\label{sec4}
\subsection{Overview of Obtained Solutions}\label{subsec4.1}
\begin{figure}[H]
\centering
\includegraphics[width=0.8\textwidth]{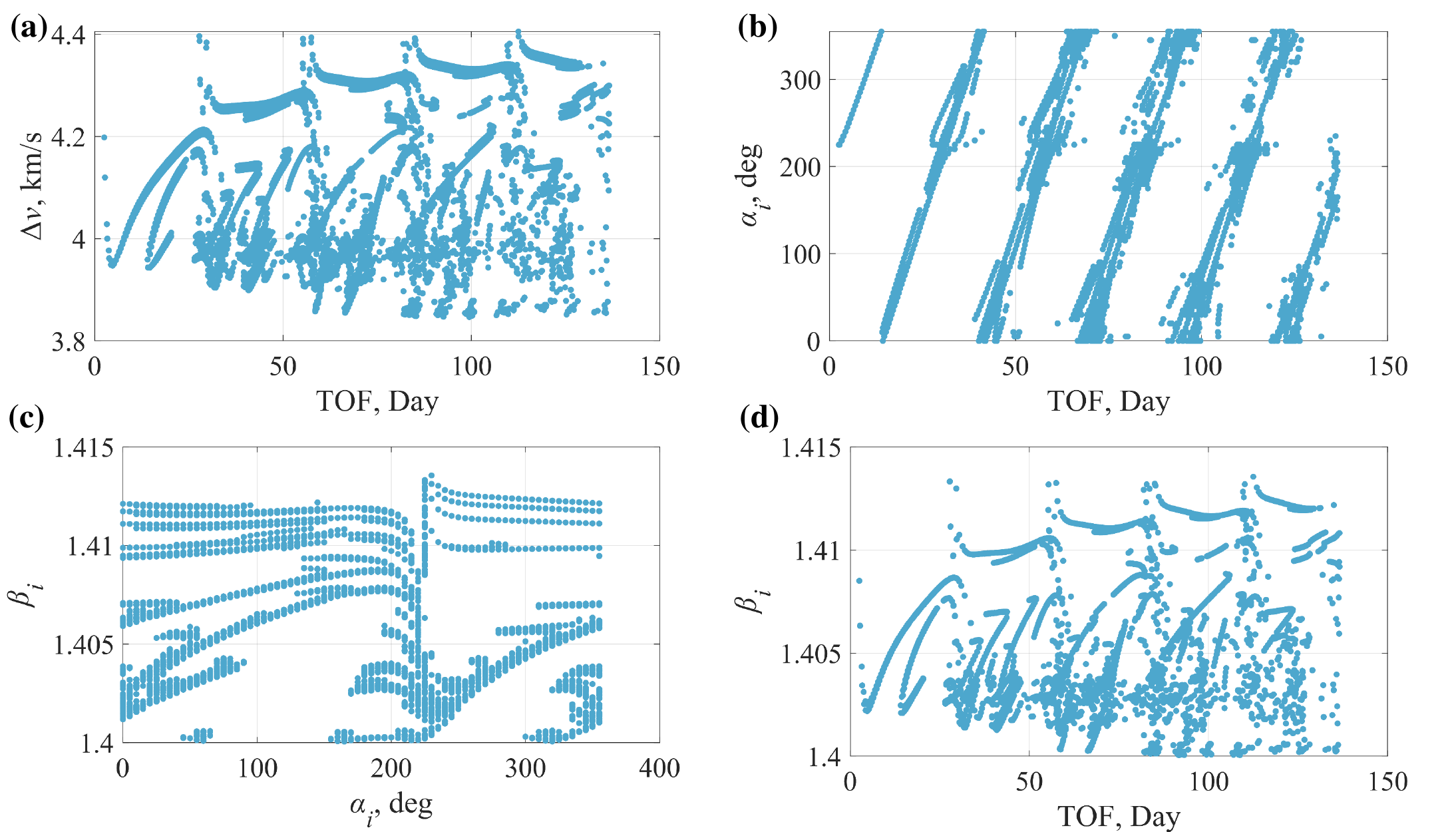}
\caption{The distribution of obtained bi-impulsive Earth-Moon transfers in the Earth-Moon PCR3BP. (a) $\left( {{\text{TOF}},{\text{ }}\Delta v} \right)$ Map; (b) $\left( {{\text{TOF}},{\text{ }}{\alpha _i}} \right)$ Map; (c) the $\left( {{\alpha _i},{\text{ }}{\beta _i}} \right)$ Map; (d) $\left( {{\text{TOF}},{\text{ }}{\beta _i}} \right)$ Map.}
\label{fig_map}
\end{figure}

Using the construction method mentioned in Section 3, we totally obtain 755654 solutions. The corresponding solution space is presented in the form of the $\left( {{\text{TOF}},{\text{ }}\Delta v} \right)$ map, the $\left( {{\text{TOF}},{\text{ }}{\alpha _i}} \right)$ map, the $\left( {{\alpha _i},{\text{ }}{\beta _i}} \right)$ map, and the $\left( {{\text{TOF}},{\text{ }}{\beta _i}} \right)$ map, as shown in Fig. \ref{fig_map}. The $\left( {{\text{TOF}},{\text{ }}\Delta v} \right)$ map of the obtained solutions is presented in Fig. \ref{fig_map} (a). The range of $\Delta v$ of the obtained solutions is $3.8466-4.4056$ km/s, and the range of TOF of the obtained solutions is $2.4975-136.5842$ Days. This type of map has been widely adopted \cite{topputo2013optimal,oshima2019low} to describe the solution space of the transfers in the multi-body problem. However, although the map provides the full information about the transfer characteristics, it does not provide more information about the distribution of the construction parameters, which can further provide useful insight into the selection of initial guesses to construct transfers. In this paper, we would like to present the distribution of the construction parameters and report an interesting, novel, or less-reported phenomenon, i.e., the discontinuous behavior of TOF distribution.

\subsection{Discontinuous Behavior of TOF distribution}\label{subsec4.2}
For the obtained solutions, the $\left( {{\text{TOF}},{\text{ }}{\alpha _i}} \right)$ map, the $\left( {{\alpha _i},{\text{ }}{\beta _i}} \right)$ map, and the $\left( {{\text{TOF}},{\text{ }}{\beta _i}} \right)$ map are presented in Fig. \ref{fig_map} (b)-(d), respectively. Differing from the distribution shown in Fig. \ref{fig_map} (c)-(d), the $\left( {{\text{TOF}},{\text{ }}{\alpha _i}} \right)$ map reveals a rather pronounced phenomenon:
\begin{enumerate}
  \item For each $\alpha_i$, the distribution of TOF reveals a discontinuous behavior. The map can be approximately divided into 6 branches, and the TOF of each branch differs by approximately one month.
  \item For each branch, $\alpha_i$ exhibits a banded distribution with respect to TOF, and the slope of $\alpha_i$ with respect to TOF is approximately equal.
\end{enumerate}

This is an interesting phenomenon that can further provide useful insight into selecting initial guesses to construct transfers, i.e., when performing the grid search, there might be no need to select the values of $\left( {{\text{TOF}},{\text{ }}{\alpha _i}} \right)$ pair located in the blank region in Fig. \ref{fig_map} (b). This phenomenon cannot be revealed by the $\left( {{\text{TOF}},{\text{ }}\Delta v} \right)$ map, as the distribution shown in Fig. \ref{fig_map} (a) reveals a continuity of TOF. How to effectively use this phenomenon to further aid in the construction of the bi-impulsive Earth-Moon transfer will be the focus of our future work.

Notably, this phenomenon is expectedly not to be limited to the specific transfer scenario (i.e., the specific setting of $h_i$ and $h_f$), as we also observed a similar phenomenon for the bi-impulsive transfer from a 36000 km circular Earth parking orbit to a 100 km circular Moon target orbit. However, the aforementioned phenomenon is revealed for bi-impulsive Earth-Moon transfers in the Earth-Moon PCR3BP. When considering the Sun-Earth/Moon PBCR4BP or other high-er-fidelity models, the distribution of $\left( {{\text{TOF}},{\text{ }}{\alpha _i}} \right)$ might be more complex due to the additional construction parameters \cite{pinelli2023neural}. How to explore the solution space of the bi-impulsive Earth-Moon transfers in these models remains an open problem.

\section{Conclusions}\label{sec5}
This paper is devoted to reporting an interesting, novel, or less-reported phenomenon about the distribution of the construction parameter of bi-impulsive Earth-Moon transfers in the Earth-Moon planar circular restricted three-body problem (PCR3BP). Firstly, the bi-impulsive transfers are constructed by grid search and trajectory correction. Then, the distribution of the obtained solutions is presented and analyzed. An interesting, novel, or less-reported phenomenon about the distribution of time of flight (TOF) and the departure phase angle at the circular Earth parking orbit has been revealed. It is observed that for each departure phase angle, the distribution of TOF reveals a discontinuous behavior, which can be divided into several branches. The TOF of each branch differs by approximately one month. This phenomenon can further pro-vide useful insight into selecting initial guesses to construct bi-impulsive transfers in the Earth-Moon PCR3BP. Future work includes the use of this phenomenon in the transfer construction in the Earth-Moon PCR3BP and the extension of this phenomenon in the higher-fidelity models.

\section*{Acknowledgements}
The authors acknowledge financial support from the National Natural Science Foundation of China (Grant Nos. 12525204, 12372044, 12302058), and the Young Elite Scientists Sponsorship Program by CAST (Grant No. 2023QNRC001).

\bibliographystyle{unsrt}

\end{document}